\documentclass[12pt]{article}
\usepackage{amssymb,latexsym,amsmath}

\begin{document}
\title{The Generalized Hodge conjecture for 1-cycles and codimension two algebraic cycles}
\author{Wenchuan Hu }






 \maketitle
 \pagestyle{myheadings}
 \markright{The GHC for 1-cycles and codimension-two algebraic cycles}

\newtheorem{Def}{Definition}[section]
\newtheorem{Th}{Theorem}[section]
\newtheorem{Prop}{Proposition}[section]
\newtheorem{Not}{Notation}[section]
\newtheorem{Lemma}{Lemma}[section]
\newtheorem{Rem}{Remark}[section]
\newtheorem{Cor}{Corollary}[section]

\def\s{\section}
\def\ss{\subsection}

\def\d{\begin{Def}}
\def\t{\begin{Th}}
\def\p{\begin{Prop}}
\def\n{\begin{Not}}
\def\la{\begin{Lemma}}
\def\r{\begin{Rem}}
\def\c{\begin{Cor}}
\def\ee{\begin{equation}}
\def\aa{\begin{eqnarray}}
\def\ya{\begin{eqnarray*}}
\def\bd{\begin{description}}

\def\ed{\end{Def}}
\def\et{\end{Th}}
\def\epo{\end{Prop}}
\def\en{\end{Not}}
\def\el{\end{Lemma}}
\def\er{\end{Rem}}
\def\ec{\end{Cor}}
\def\eee{\end{equation}}
\def\eaa{\end{eqnarray}}
\def\ey{\end{eqnarray*}}
\def\ebd{\end{description}}

\def\nn{\nonumber}
\def\bp{{\bf Proof.}\hspace{2mm}}
\def\qe{\hfill$\Box$}
\def\lj{\langle}
\def\rj{\rangle}
\def\dd{\diamond}
\def\ox{\mbox{}}
\def\lb{\label}
\def\rel{\;{\rm rel.}\;}
\def\vp{\varepsilon}
\def\ep{\epsilon}
\def\mod{\;{\rm mod}\;}
\def\exp{{\rm exp}\;}
\def\Lie{{\rm Lie}}
\def\dim{{\rm dim}}
\def\im{{\rm im}\;}
\def\Lag{{\rm Lag}}
\def\Gr{{\rm Gr}}
\def\span{{\rm span}}
\def\Spin{{\rm Spin}}
\def\sign{{\rm sign}\;}
\def\Supp{{\rm Supp}\;}
\def\Sp{{\rm Sp}\;}
\def\ind{{\rm ind}\;}
\def\rank{{\rm rank}\;}
\def\Sg{{\Sp(2n,\C)}}
\def\Na{{ \mathcal{N}}}
\def\det{{\rm det}\;}
\def\dist{{\rm dist}}
\def\deg{{\rm deg}}
\def\Tr{{\rm Tr}\;}
\def\ker{{\rm ker}\;}
\def\Vect{{\rm Vect}}
\def\H{{\bf H}}
\def\K{{\rm K}}
\def\R{{\mathbb{R}}}
\def\C{{\mathbb{C}}}
\def\Z{{\mathbb{Z}}}
\def\N{{\bf N}}
\def\F{{\bf F}}
\def\Da{{\bf D}}
\def\A{{\bf A}}
\def\La{{\bf L}}
\def\x{{\bf x}}
\def\y{{\bf y}}
\def\Ga{{\cal G}}
\def\Ha{{\cal H}}
\def\L{{\cal L}}
\def\Pa{{\cal P}}
\def\Ua{{\cal U}}
\def\E{{\rm E}}
\def\J{{\mathcal{J}}}

\def\m{{\rm m}}
\def\ch{{\rm ch}}
\def\gl{{\rm gl}}
\def\Gl{{\rm Gl}}
\def\Sp{{\rm Sp}}
\def\sf{{\rm sf}}
\def\U{{\rm U}}
\def\O{{\rm O}}
\def\F{{\rm F}}
\def\P{{\rm P}}
\def\D{{\rm D}}
\def\T{{\rm T}}
\def\Sa{{\rm S}}

\begin{abstract} In this paper, we prove that the statement:
``The (Generalized) Hodge
Conjecture holds for codimension-two cycles on a smooth
projective variety $X$" is a birationally invariant statement,
that is, if the statement is true for $X$, it is also true for all
smooth  varieties $X'$ which are birationally equivalent to $X$. We also prove the analogous
result for 1-cycles. As direct corollaries, the Hodge Conjecture holds for
smooth rational projective manifolds with dimension less than or equal to
five, and, the Generalized Hodge Conjecture holds for smooth rational projective
manifolds with dimension less than or equal to four.
\end{abstract}

\begin{center}{\bf \tableofcontents}\end {center}

\medskip
\s {Introduction} In this paper, all varieties are defined over
$\mathbb{C}$. Let $X$ be a smooth projective variety with dimension
$n$. Let ${\mathcal Z}_p(X)$ be the space of algebraic $p$-cycles on
$X$. Set ${\mathcal Z}^{n-p}(X)\equiv{\mathcal Z}_p(X)$. There is a
natural map

$$cl_q:{\mathcal Z}^q(X)\rightarrow H^{2q}(X,\Z)$$
called \textbf{the cycle class map}.

\medskip
Tensoring with ${\mathbb{Q}}$, we have

$$cl_q\otimes{\mathbb{Q}}:{\mathcal Z}^q(X)\otimes{\mathbb{Q}}\rightarrow H^{2q}(X,\mathbb{Q}).$$

It is well known that  $cl_q({\mathcal Z}_q(X))\subseteq
H^{q,q}(X)\cap \rho(H^{2q}(X,\Z))$, where $\rho:
H^{2q}(X,\Z)\rightarrow H^{2q}(X,\mathbb{C})$  is the coefficient
homomorphism and $H^{q,q}(X)$ denotes the $(q,q)$-component in the
Hodge decomposition (cf.\cite{Griffiths-Harris}, \cite{Lewis1}).
There are known examples where $cl_q({\mathcal Z}_q(X))\neq
H^{q,q}(X)\cap \rho (H^{2q}(X,\Z))$ (cf.
\cite{Ballico-Catanese-Ciliberto} p.134-125], \cite{Lewis2}). We
recall:

\medskip
\medskip
\noindent {\bf The Hodge Conjecture} (for codimension-q cycles):
{\sl The rational cycle class map

$$cl_q\otimes{\mathbb{Q}}:{\mathcal
Z}^q(X)\otimes{\mathbb{Q}}\rightarrow H^{q,q}(X)\cap
H^{2q}(X,\mathbb{Q})$$ is surjective. }

\medskip
\medskip
\noindent {\bf The Hodge Conjecture over $\Z$}: {\sl The rational
cycle class map

$$cl_q:{\mathcal Z}^q(X)\rightarrow H^{q,q}(X)\cap
\rho (H^{2q}(X,\Z))$$ is surjective. }

\medskip
We shall denote by $\rm Hodge^{q,q}(X,\mathbb{Q})$ the statement
that: `` The Hodge Conjecture for codimension-q cycles is true for
$X$". Similarly, we denote by $\rm Hodge^{q,q}(X,\bf Z)$ the
corresponding statement for the Hodge Conjecture over $\Z$.

\medskip
More generally, we can define a  filtration on $H_k(X,{\mathbb{Q}})$
as follows:

{\Def {\rm (\cite{Friedlander-Mazur},\S 7])} Denote by
$\tilde{F}_pH_k(X,{\mathbb{Q}})\subseteq H_k(X,{\mathbb{Q}})$  the
maximal sub-(Mixed) Hodge structure of span $k-2p$. (See
\cite{Grothendieck} and \cite{Friedlander-Mazur}.) The
sub-${\mathbb{Q}}$ vector spaces $\tilde{F}_pH_k(X,{\mathbb{Q}})$
form a decreasing filtration of sub-Hodge structures:

$$\cdots\subseteq \tilde{F}_pH_k(X,{\mathbb{Q}})\subseteq
\tilde{F}_{p-1}H_k(X,{\mathbb{Q}}) \subseteq\cdots\subseteq
\tilde{F}_0H_k(X,{\mathbb{Q}})\subseteq H_k(X,{\mathbb{Q}})$$ and
$\tilde{F}_pH_k(X,{\mathbb{Q}})$ vanishes if $2p>k$. This filtration
is called the \textbf{Hodge filtration}.}

\medskip
A homological version of the arithmetic filtration (see
[\cite{Lewis1},\S 7]) is given in the following definition:

{\Def {\rm (\cite{Friedlander-Mazur},\S 7])} Denote by
$G_pH_k(X,{\mathbb{Q}})\subseteq H_k(X,{\mathbb{Q}})$ the
$\mathbb{Q}$-vector subspace of $H_k(X,{\mathbb{Q}})$ generated by
the images of mappings $H_k(Y,{\mathbb{Q}})\rightarrow
H_k(X,{\mathbb{Q}})$,induced from all morphisms $Y\rightarrow X$ of
varieties of dimension$\leq k-p$. The subspaces
$G_pH_k(X,{\mathbb{Q}})$ also form a decreasing filtration called
the \textbf{geometric filtration}:

$$\cdots\subseteq
G_pH_k(X,{\mathbb{Q}})\subseteq G_{p-1}H_k(X,{\mathbb{Q}})
\subseteq\cdots\subseteq G_0H_k(X,{\mathbb{Q}})\subseteq H_k(X,{
\mathbb{Q}}).$$ }

Since $X$ is smooth, the Weak Lefschetz Theorem implies that
$G_0H_k(X,{\mathbb{Q}})= H_k(X,{\mathbb{Q}})$. Since
$H_k(Y,{\mathbb{Q}})$ vanishes for $k$ greater than twice the
dimension of $Y$, $G_pH_k(X,{\mathbb{Q}})$ vanishes if $2p>k$.

It was proved in \cite{Grothendieck} that, for any smooth variety
$X$, the geometric filtration is finer than the Hodge filtration,
i.e.,
$G_pH_k(X,{\mathbb{Q}})\subseteq\tilde{F}_pH_k(X,{\mathbb{Q}})$, for
all $p$ and $k$.

\medskip
\noindent  \textbf{The Generalized Hodge Conjecture}: For any smooth
variety $X$,

\begin{equation} \label{eq:GHC}
G_pH_k(X,{\mathbb{Q}})=\tilde{F}_pH_k(X,{\mathbb{Q}})
\end{equation}

\medskip
\noindent for all $p$ and $k$. Using the notation given in
\cite{Lewis1}, we denote by $\widetilde{GHC}(p,k,X)$ the assertion
that (\ref{eq:GHC}) is true.

\medskip
{\Def The  \textbf{Lawson homology} $L_pH_k(X)$ of $p$-cycles is
defined by

$$L_pH_k(X) = \pi_{k-2p}({\mathcal Z}_p(X)) \quad for\quad
k\geq 2p\geq 0,
$$

\medskip
\noindent where ${\mathcal Z}_p(X)$ is provided with a natural
topology (cf. \cite{Friedlander}, \cite{Lawson1}).} For general
background, the reader is referred to Lawson' survey paper
\cite{Lawson2}.

\medskip
There are two special cases.

\begin{enumerate}
\item[(a)] If $p=0$, then for all $k\geq 0$,
$L_0H_k(X)\cong H_k(X,\Z)$ by Dold-Thom Theorem \cite{Dold-Thom}.
\item[(b)] If $k=2p$, then
$L_pH_{2p}(X)={\mathcal Z}_p(X)/{{\mathcal Z}_p(X)}_{alg}$, where
${{\mathcal Z}_p(X)}_{alg}$ denotes the algebraic $p$-cycles on $X$
which are algebraic equivalent to zero.
\end{enumerate}

\medskip
In \cite{Friedlander-Mazur}, Friedlander and Mazur showed that there
are natural maps, called \textbf{cycle class maps}

$$ \Phi_{p,k}:L_pH_{k}(X)\rightarrow H_{k}(X,\Z). $$
Define $$L_pH_{k}(X)_{hom}:={\rm
ker}\{\Phi_{p,k}:L_pH_{k}(X)\rightarrow H_{k}(X,\Z)\}.$$

$$T_pH_k(X):={\rm Im}
\{\Phi_{p,k}:L_pH_k(X)\rightarrow H_k(X,\Z)\}$$ and
$$T_pH_k(X,{\mathbb{Q}}):=T_pH_k(X)\otimes {\mathbb{Q}}.$$

\medskip
It was proved in [\cite{Friedlander-Mazur},\S7] that, for any smooth
variety $X$, $T_pH_k(X, {\mathbb{Q}})\subseteq
G_pH_k(X,{\mathbb{Q}})$ for all $p$ and $k$. Hence

\begin{equation} \label{eq:TopGeoHodge}
T_pH_k(X, {\mathbb{Q}})\subseteq G_pH_k(X,{\mathbb{Q}})\subseteq
\tilde{F}_pH_k(X, {\mathbb{Q}}).
\end{equation}

\medskip
In this paper, we will use the tools in Lawson homology  and the
methods given in \cite{author1} to show the following main result:

{\Th Let $X$ be a smooth projective variety. If the Hodge conjecture
for codimension 2 cycles over $\Z$ holds for $X$, i.e., if we have
$\rm Hodge^{2,2}(X,\Z)$, then it holds for any smooth projective
variety $X^{\prime}$ birational to $X$. That is, $\rm
Hodge^{2,2}(X,\Z)$ is a birationally invariant assertion for smooth
varieties $X$.}

{\Rem The above theorem remains true if  $\Z$ is replaced by
$\mathbb{Q}$. Since $\rm Hodge^{2,2}(X,\mathbb{Q})$ implies $\rm
Hodge^{n-2,n-2}(X,\mathbb{Q})$ for $n\geq 4$ (cf. [\cite{Lewis1},
p.91]), $\rm Hodge^{n-2,n-2}(X,\mathbb{Q})$ is also a birationally
invariant property of smooth n-dimensional varieties $X$. }

\medskip
As a corollary, we have

{\Cor If $X$ is a rational manifold with ${\rm dim}(X)\leq 5$, then
the Hodge conjecture $\rm Hodge^{p,p}(X,\mathbb{Q})$ is true for
$1\leq p \leq {\rm dim}(X)$. In fact, $\rm
Hodge^{p,p}(X^{\prime},\bf Z)$ is true except possibly for $p=3,
\dim(X)=5$.}

{\Rem By using the technique of the  diagonal decomposition, Bloch
and Srinivas \cite{Bloch-Srinivas} showed that, for any smooth
projective variety $X$, $\rm Hodge^{2,2}(X,\mathbb{Q})$ holds if the
Chow group of 0-cycles ${\rm Ch}_0(X)\cong \Z$ . Laterveer
\cite{Laterveer} generalized this technique and showed the Hodge
Conjecture holds for a class of projective manifolds with small chow
groups. }

\medskip
{\Cor Let $X$ be a smooth projective variety of dimension $\leq 5$
such that the Hodge Conjecture is known to be true, i.e., $\rm
Hodge^{p,p}(X,\mathbb{Q})$ holds for all $p$. Then the Hodge
Conjecture holds for all smooth projective varieties $X'$ which are
birationally equivalent to $X$. Non-rational examples of such an $X$
include general abelian varieties or the product of at most five
elliptic curves. For more examples, the reader is referred to the
survey book \cite{Lewis1}. }

\medskip
Our second main result is the following

{\Th The assertion $\widetilde{GHC}(n-2,k,X)$ is a birationally
invariant property of smooth $n$-dimensional varieties $X$ when
$k\geq 2(n-2)$. More precisely, if $\widetilde{GHC}(n-2,k,X)$ holds
for a smooth variety $X$, then $\widetilde{GHC}(n-2,k,X')$ holds for
any smooth variety $X'$ birational to $X$.}

\medskip
We also show that

{\Prop The assertion that
``$T_{n-2}H_k(X,{\mathbb{Q}})=\tilde{F}_{n-2}H_k(X,{\mathbb{Q}})$
holds" is a birationally invariant property of smooth
$n$-dimensional varieties $X$ when $k\geq 2(n-2)$.}

\medskip
Similarly, for 1-cycles,  we can show the following.

{\Prop  For integer $k\geq 2$, the assertion that
``$T_1H_k(X,{\mathbb{Q}})=\tilde{F}_1H_k(X,{\mathbb{Q}})$ holds" is
a birationally invariant property of smooth $n$-dimensional
varieties $X$.}

\medskip
and

{\Th For any integer $k\geq 2$, the assertion $\widetilde{GHC}(1,k,X)$  is a
birationally invariant property of smooth varieties $X$.}

{\Rem For the case $k=\dim(X)$, Lewis has already obtained this
result in \cite{Lewis1}.}

{\Cor For any smooth rational variety $X$ with $\dim (X)\leq 4$, the
Generalized Hodge Conjecture holds.}

\medskip
The main tools used to prove this result are: the long exact
localization sequence given by Lima-Filho in \cite{Lima-Filho}, the
explicit formula for Lawson homology of codimension-one cycles on a
smooth projective manifold given by Friedlander in
\cite{Friedlander}, and the weak factorization theorem proved by
Wlodarczyk in \cite{Wlodarczyk} and in \cite{AKMW}.

\medskip
\s {The Proof of the Main Theorems }  Let $X$ be a smooth projective
manifold of dimension $n$. In the following, we will denote by
$H_{p,q}(X)$ the image of $H^{n-p,n-q}(X)$ under the Poincare
duality isomorphism $H^{2n-p-q}(X,\mathbb{C})\cong
H_{p+q}(X,\mathbb{C})$.

\medskip
Let $X$ be a smooth projective manifold and $i_0:Y\hookrightarrow X$
be a smooth subvariety of codimension $r$. Let
$\sigma:\tilde{X}_Y\rightarrow X$ be the blowup of $X$ along $Y$,
$i:D:=\sigma^{-1}(Y)\hookrightarrow \tilde{X}_Y$ the exceptional
divisor of the blowing up, and $\pi:D\rightarrow Y$ the restriction
of $\sigma$ to $D$. Set $U:= X-Y\cong \tilde{X}_Y - D$. Denote
by $j_0$ the inclusion $U\subset X$ and $j$ the inclusion $U\subset
\tilde{X}_Y$.

\medskip
Now I list the Lemmas and Corollaries given in \cite{author1}.

{\Lemma \label{eq:bpdiaglh}For each $p$, we have the following
commutative diagram
$$
\begin{array}{ccccccccccc}
\cdots\rightarrow & L_pH_k(D) & \stackrel{i_*}{\rightarrow} &
L_pH_k({\tilde{X}_Y})
 & \stackrel{j^*}{\rightarrow} & L_pH_k(U) & \stackrel{\delta_*}{\rightarrow}
 & L_pH_{k-1}(D) & \rightarrow & \cdots & \\

 & \downarrow \pi_*&   & \downarrow \sigma_*&   & \downarrow  \cong&   & \downarrow \pi_*&   &  &\\

 \cdots\rightarrow & L_pH_k(Y) & \stackrel{(i_0)*}{\rightarrow} &
L_pH_k({X}) & \stackrel{j_0^*}{\rightarrow}&
L_pH_k(U)&\stackrel{(\delta_0)_*}{\rightarrow} & L_pH_{k-1}(Y) &
\rightarrow & \cdots&

\end{array}
$$
}

{\Rem Since $\pi_* $ is surjective (this follows from the explicit
formula for the Lawson homology of $D$, i.e., the Projective Bundle
Theorem in \cite{Friedlander-Gabber}), it is easy to see that
$\sigma_*$ is surjective.}

{\Cor \label{sec:bpdiaghom}If $p=0$, then we have the following
commutative diagram

$$
\begin{array}{ccccccccccc}
\cdots\rightarrow & H_k(D) & \stackrel{i_*}{\rightarrow} &
H_k({\tilde{X}_Y}) & \stackrel{j^*}{\rightarrow} & H_k^{BM}(U) &
\stackrel{\delta_*}{\rightarrow} & H_{k-1}(D) & \rightarrow & \cdots
& \\

& \downarrow \pi_*&   & \downarrow \sigma_*&   & \downarrow  \cong&
& \downarrow \pi_*&   &  &\\

\cdots\rightarrow & H_k(Y) & \stackrel{(i_0)*}{\rightarrow} &
H_k({X}) & \stackrel{j_0^*}{\rightarrow} &
H_k^{BM}(U)&\stackrel{(\delta_0)_*}{\rightarrow} & H_{k-1}(Y) &
\rightarrow & \cdots&

\end{array}
$$
Moreover, if $x\in H_k(D)$ vanishes under $\pi_*$ and $i_*$, then
$x=0\in H_k(D)$. }

{\Cor If $p=n-2$, then we have the following commutative diagram
$$
\begin{array}{ccccccccccc}
\cdots\rightarrow & L_{n-2}H_k(D) & \stackrel{i_*}{\rightarrow} &
L_{n-2}H_k({\tilde{X}_Y}) & \stackrel{j^*}{\rightarrow} &
L_{n-2}H_k(U) & \stackrel{\delta_*}{\rightarrow} & L_{n-2}H_{k-1}(D)
& \rightarrow & \cdots & \\

& \downarrow \pi_*&   & \downarrow \sigma_*&   & \downarrow  \cong&
& \downarrow \pi_*&   &  &\\

\cdots\rightarrow & L_{n-2}H_k(Y) & \stackrel{(i_0)*}{\rightarrow} &
L_{n-2}H_k({X}) & \stackrel{j_0^*}{\rightarrow} &
L_{n-2}H_k(U)&\stackrel{(\delta_0)_*}{\rightarrow} &
L_{n-2}H_{k-1}(Y) & \rightarrow & \cdots&

\end{array}
$$
}

{\Lemma  For each $p$, we have the following commutative diagram
$$
\begin{array}{ccccccccccc}
\cdots\rightarrow & L_pH_k(D) & \stackrel{i_*}{\rightarrow} &
L_pH_k({\tilde{X}_Y}) & \stackrel{j^*}{\rightarrow} & L_pH_k(U) &
\stackrel{\delta_*}{\rightarrow} & L_pH_{k-1}(D) &
\rightarrow & \cdots & \\

& \downarrow \Phi_{p,k}&   & \downarrow \Phi_{p,k}&   & \downarrow
\Phi_{p,k}&   & \downarrow \Phi_{p,k-1}&   &  &\\

\cdots\rightarrow & H_k(D) & \stackrel{i_*}{\rightarrow} &
H_k({\tilde{X}_Y}) & \stackrel{j^*}{\rightarrow} &
H_k^{BM}(U)&\stackrel{\delta_*}{\rightarrow} & H_{k-1}(D) &
\rightarrow & \cdots&

\end{array}
$$
In particular, it is true for $p=1,n-2$.}

{\Lemma For each $p$, we have the following commutative diagram

$$
\begin{array}{ccccccccccc}
\cdots\rightarrow & L_pH_k(Y) & \stackrel{(i_0)_*}{\rightarrow} &
L_pH_k({{X}}) & \stackrel{j^*}{\rightarrow}& L_pH_k(U) &
\stackrel{(\delta_0)_*}{\rightarrow} & L_pH_{k-1}(Y) &
\rightarrow & \cdots & \\

& \downarrow \Phi_{p,k}&   & \downarrow \Phi_{p,k}&   & \downarrow
\Phi_{p,k}&   & \downarrow \Phi_{p,k-1}&   &  &\\

\cdots\rightarrow & H_k(Y) & \stackrel{(i_0)_*}{\rightarrow} &
H_k({X}) & \stackrel{j^*}{\rightarrow} &
H_k^{BM}(U)&\stackrel{(\delta_0)_*}{\rightarrow} & H_{k-1}(Y) &
\rightarrow & \cdots&

\end{array}
$$
In particular, it is true for $p=1,n-2$.}

{\Rem All the commutative diagrams of long exact sequences remain
commutative and exact after tensoring with ${\mathbb{Q}}$. We will
use these Lemmas and corollaries with rational coefficients. }

\medskip
The following result  proved by Friedlander will be used several
times:

{\Th {\rm (Friedlander \cite{Friedlander})} \label{sec:codim2} Let
$X$ be any smooth projective variety of dimension $n$. Then we have
the following isomorphisms

$$
\left\{
\begin{array}{l}
 L_{n-1}H_{2n}(X)\cong \Z,\\
 L_{n-1}H_{2n-1}(X)\cong H_{2n-1}(X,\Z),\\
 L_{n-1}H_{2n-2}(X)\cong H_{n-1,n-1}(X,\Z)=NS(X)\\
L_{n-1}H_{k}(X)=0 \quad for\quad k> 2n.\\

\end{array}
\right.
$$
where $NS(X)$ is the N\'{e}ron-Severi group of $X$. }


\ss{The Proof of Theorem 1.1 for a Blowup}  In what follows we drop
reference to the coefficient homomorphism $\rho$, and denote by
$H_k(X,\Z)$ its image in $H_k(X,\mathbb{C})$.

\medskip
There are two cases to consider:

\medskip
\noindent\textbf{Case 1:} If $\Phi_{n-2,2(n-2)}:
L_{n-2}H_{2(n-2)}(X)\rightarrow H_{2(n-2)}(X,\Z)\cap H_{n-2,n-2}(X)
$ is surjective, we will show that $\Phi_{n-2,2(n-2)}:
L_{n-2}H_{2(n-2)}(\tilde{X}_Y)\rightarrow
H_{2(n-2)}(\tilde{X}_Y,\Z)\cap H_{n-2,n-2}(\tilde{X}_Y)$ is also
surjective.

\medskip
Let $b\in H_{n-2,n-2}(\tilde{X}_Y)\cap H_{2(n-2)}(\tilde{X}_Y,\Z)$.
Set $a\equiv \sigma_*(b)\in H_{2(n-2)}(X,\Z)$. Since $\sigma_*$
preserves the type, we have $a\in H_{2(n-2)}(X,\Z)\cap
H_{n-2,n-2}(X)$. Now by assumption, there exists an element
$\tilde{a}\in L_{n-2}H_{2(n-2)}(X)$ such that
$\Phi_{n-2,2(n-2)}(\tilde{a})=a$. Now since
$L_{n-2}H_{2(n-2)}(\tilde{X}_Y)\rightarrow L_{n-2}H_{2(n-2)}(X)$ is
surjective, there exists an element $\tilde{b}\in
L_{n-2}H_{2(n-2)}(\tilde{X}_Y)$ such that
$\sigma_*(\tilde{b})=\tilde{a}$. Now
$\Phi_{n-2,2(n-2)}(\tilde{b})-b$ is mapped to zero under $\sigma_*$
on $H_{2(n-2)}(\tilde{X}_Y,\Z)$. By the commutative diagram in the
long exact sequences in Corollary \ref{sec:bpdiaghom}, there exists
an element $c\in H_{2(n-2)}(D,\Z)$ such that
$i_*(c)=\Phi_{n-2,2(n-2)}(\tilde{b})-b$. Using Corollary
\ref{sec:bpdiaghom} once again, we have $\pi_*(c)=0$. This follows
from the fact that $\dim (Y)=n-r\leq n-2$ and hence
$(i_0)_*:H_{2(n-2)}(Y,\Z)\rightarrow H_{2(n-2)}(X,\Z)$ is injective.
From the blowup formula for the singular homology, $i_*|_{\ker
\pi_*}$ is injective. Now by assumption, $b$ and $\tilde{b}$ are
non-torsion elements. Hence $c$ is not a torsion element in
$H_{2(n-2)}(D,\Z)$, i.e.,$c\in H_{2(n-2)}(D,\Z)_{\rm free}$, the
torsion free part of $H_{2(n-2)}(D,\Z)$.

Since $i_*$ preserves the type, we have the following

\medskip
{\bf Claim}: {\sl $c\in H_{2(n-2)}(D,\Z)\cap H_{n-2,n-2}(D)$.}

\medskip
\bp Since $H_{2(n-2)}(D,\Z)_{\rm free}\subset
H_{2(n-2)}(D,{\mathbb{C}}) =H_{n-2,n-2}(D)\oplus
H_{n-1,n-3}(D)\oplus H_{n-3,n-1}(D)$. Now $c=c_0+c_1+\bar{c_1}\in
H_{2(n-2)}(D,{\mathbb{C}})$ such that $c_0\in H_{n-2,n-2}(D)$,
$c_1\in H_{n-1,n-3}(D)$ and hence $\bar{c_1}\in H_{n-3,n-1}(D)$.
Note that the complexification of $i_*$ is the map $i_*\otimes
\mathbb{C}: H_{2(n-2)}(D,{\mathbb{C}})\rightarrow
H_{2(n-2)}(\tilde{X}_Y,{\mathbb{C}})$. If $i_*\otimes \C (c_1)=0$,
we have $c_1=0$. In fact, $i_*\otimes \C (c_1)=0$ and the exactness
of the long exact sequence in the upper row in Corollary
\ref{sec:bpdiaghom} implies that an element $d\in
H^{BM}_{2(n-2)+1}(U,{\mathbb{C}})$ such that $\delta_*(d)=c_1$. We
use the commutative diagram in Corollary \ref{sec:bpdiaghom} again.
From the commutativity of the diagram in Corollary
\ref{sec:bpdiaghom}, we have the image of $d$ under the boundary map
$(\delta_0)_*$ must zero in $H_{2(n-2)}(Y,{\mathbb{C}})$.  This
follows from the fact that the complex dimension of $\dim (Y)\leq
n-2$ and the Hodge type of $d$ is of type $(n-1,n-3)$. Now by the
exactness of the long exact sequence in the lower row in Corollary
\ref{sec:bpdiaghom}, there exists an element $e\in
H_{2(n-2)+1}(X,{\mathbb{C}})$ such that $j_0^*(e)=d$. It is
well-known that
$\sigma_*:H_{2(n-2)+1}(\tilde{X}_Y,{\mathbb{C}})\rightarrow
H_{2(n-2)+1}(X,{\mathbb{C}})$ is surjective. Therefore, there exists
$\tilde{e}\in H_{2(n-2)+1}(\tilde{X}_Y,{\mathbb{C}})$ such that
$\sigma_*(\tilde{e})=e$. We get $d=j^*(\tilde{e})$ and hence
$c_1=0\in H_{2(n-2)}(D,{\mathbb{C}})$ by the exactness of the the
upper row sequence in  Corollary \ref{sec:bpdiaghom}. This implies
$\bar{c_1}=0$ and hence $c\in H_{n-2,n-2}(D)$. This finishes the
proof of the claim.

\qe

\medskip
Since $\dim D=n-1$, hence by Theorem \ref{sec:codim2}, the map
$\Phi_{n-2,2(n-2)}: L_{n-2}H_{2(n-2)}(D)\rightarrow
H_{2(n-2)}(D,\Z)\cap H_{n-2,n-2}(D) $ is an isomorphism. Set
$\tilde{c}\equiv\Phi_{n-2,2(n-2)}(c)$. Therefore,
$\Phi_{n-2,2(n-2)}\{\tilde{b}-i_*(\tilde{c})\}=b$. Hence
$\Phi_{n-2,2(n-2)}:L_{n-2}H_{2(n-2)}(\tilde{X}_Y)\rightarrow
H_{2(n-2)}(\tilde{X}_Y,\Z)\cap H_{n-2,n-2}(\tilde{X}_Y)$ is
surjective .

\medskip
On the other hand, we need to show

\medskip

\noindent\textbf{Case 2:} If $\Phi_{n-2,2(n-2)}:
L_{n-2}H_{2(n-2)}(\tilde{X}_Y)\rightarrow
H_{2(n-2)}(\tilde{X}_Y,\Z)\cap H_{n-2,n-2}(\tilde{X}_Y)$ is
surjective, then $\Phi_{n-2,2(n-2)}: L_{n-2}H_{2(n-2)}(X)\rightarrow
H_{2(n-2)}(X,\Z)\cap H_{n-2,n-2}(X)$ is also surjective.

\medskip
This part is relatively easy. Let $a\in H_{2(n-2)}(X)\cap
H_{n-2,n-2}(X)$. Since $\sigma_*:H_{2(n-2)}(\tilde{X}_Y,\Z)
\rightarrow H_{2(n-2)}(X,\Z)$ is surjective and $\sigma_*\otimes
\mathbb{C} :H_{2(n-2)}(\tilde{X}_Y,\mathbb{C}) \rightarrow
H_{2(n-2)}(X,\C)$ preserves the Hodge type, there exists an element
$b\in H_{2(n-2)}(\tilde{X}_Y,\Z)\cap H_{n-2,n-2}(\tilde{X}_Y)$ such
that $\sigma_*(b)=a$. Now by assumption, we have an element
$\tilde{b}\in L_{n-2}H_{2(n-2)}(\tilde{X}_Y)$ such that
$\Phi_{n-2,2(n-2)}(\tilde{b})=b$. Set $\tilde{a}\equiv
\sigma_*(\tilde{b})$. Then from the commutative of the diagram, we
have $\Phi_{n-2,2(n-2)}(\tilde{a})=a$. This is exactly the
surjectivity in this case.

This completes the proof for a blowup along a smooth codimension at
least two subvariety $Y$ in $X$.

\qe


\ss{The Proof of Theorem 1.2 for a Blowup}  Now we have the
following:

{\Prop \label{sec:tophdgcod2}The assertion that
``$T_{n-2}H_k(X,{\mathbb{Q}})=\tilde{F}_{n-2}H_k(X,{\mathbb{Q}})$
holds" is a birationally invariant property of smooth
$n$-dimensional varieties $X$ when $k\geq 2(n-2)$.}

\medskip
\bp There are two cases to consider:

\medskip
\textbf{Case A:} If $\Phi_{n-2,k}: L_{n-2}H_{k}(X)\otimes {
\mathbb{Q}}\rightarrow \tilde{F}_{n-2}H_{k}(X,{\mathbb{Q}}) $ is
surjective, we want to show $\Phi_{n-2,k}:
L_{n-2}H_{k}(\tilde{X}_Y)\otimes {\mathbb{Q}}\rightarrow
{\tilde{F}}_{n-2}H_{k}(\tilde{X}_Y,{\mathbb{Q}}) $ is also
surjective.

\medskip
Let $a\in \tilde{F}_{n-2}H_{k}(\tilde{X}_Y,{\mathbb{Q}})$, set
$b=\sigma_*(a)\in \tilde{F}_{n-2}H_{k}(X,{\mathbb{Q}})$. By
assumption, there exists $\tilde{b}\in L_{n-2}H_{k}(X,{\mathbb{Q}})$
such that $\Phi_{n-2,k}(\tilde{b})=b$. By the blowup formula in
Lawson homology (see \cite{author1}), we know that $\sigma_*:
L_{n-2}H_{k}(\tilde{X}_Y,{\mathbb{Q}})\rightarrow
L_{n-2}H_{k}(X,{\mathbb{Q}})$ is surjective, there exists an element
$\tilde{a}\in L_{n-2}H_{k}(\tilde{X}_Y,{\mathbb{Q}})$ such that
$\sigma_*(\tilde{a})=\tilde{b}$. By the commutative diagram in Lemma
\ref{eq:bpdiaglh} and Corollary \ref{sec:bpdiaghom}, we have
$j^*(\Phi_{n-2,k}(\tilde{a})-a)=0\in H_k^{BM}(U,{\mathbb{Q}})$. The
exactness of the localization sequence in the rows in Corollary
\ref{sec:bpdiaghom} implies that there exists an element $c\in
H_k(D, {\mathbb{Q}})$ such that $i_*(c)=\Phi_{n-2,k}(\tilde{a})-a$.
Since the $\dim(D)=n-1$ and $D$ is smooth, by Theorem
\ref{sec:codim2}, we know the natural transformation $\Phi_{n-2,k}:
L_{n-2}H_k(D)\rightarrow H_k(D)$ is an isomorphism for $k\geq
2(n-2)+1$. Hence $\Phi_{n-2,k}: L_{n-2}H_k(D)\otimes {
\mathbb{Q}}\cong H_k(D,{\mathbb{Q}})$. Therefore there exists
$\tilde{c}\in L_{n-2}H_k(D)\otimes {\mathbb{Q}}$ such that
$\Phi_{n-2,k}(\tilde{c})=c$. Now it is obvious that
$\Phi_{n-2,k}(\tilde{a}-i_*(\tilde{c}))=a$. The proof of the case
$k=2(n-p)$ is from the proof of Theorem 1.1. This is the
surjectivity as we want.

\medskip
\textbf{Case B:} If $\Phi_{n-2,k}: L_{n-2}H_{k}(\tilde{X}_Y)\otimes
{\mathbb{Q}}\rightarrow \tilde{F}_{n-2}H_{k}(\tilde{X}_Y,
{\mathbb{Q}}) $ is surjective, we want to show $\Phi_{n-2,k}:
L_{n-2}H_{k}(X)\otimes {\mathbb{Q}}) \rightarrow
\tilde{F}_{n-2}H_{k}(X,{\mathbb{Q}}) $ is also surjective. We can
use an argument similar to the \textbf{Case 2} above. Suppose $b\in
\tilde{F}_{n-2}H_{k}(X,{\mathbb{Q}})$. Then there exists a
$\tilde{b}\in \tilde{F}_{n-2}H_{k}(\tilde{X}_Y,{ \mathbb{Q}})$ such
that $\sigma_*(\tilde{b})=b$ by the blowup formula for the singular
homology with ${\mathbb{Q}}$-coefficients. By assumption, there
exists an $\tilde{a}\in L_{n-2}H_{k}(\tilde{X}_Y)\otimes
{\mathbb{Q}}$ such that $\Phi_{n-2,k}(\tilde{a})=\tilde{b}$. Set
$a=\sigma_*(\tilde{a})$. Then $a\in L_{n-2}H_{k}({X})\otimes
{\mathbb{Q}}$ and $\Phi_{n-2,k}(a)=b$. This finishes the proof of
the surjectivity in this case.

\qe

\medskip
Now we give the proof of Theorem 1.2. First, we suppose that
$G_{n-2}H_k(X,{\mathbb{Q}})=\tilde{F}_{n-2}H_k(X,{\mathbb{Q}})$. We
will show
$G_{n-2}H_k(\tilde{X}_Y,{\mathbb{Q}})=\tilde{F}_{n-2}H_k(\tilde{X}_Y,{\mathbb{Q}})$
case by case.

For $k>2n$, $ \tilde{F}_{n-2}H_{k}(\tilde{X}_Y)=0$ and hence nothing
needs to be proved.

For $k=2n$,
$G_{n-2}H_{k}(\tilde{X}_Y)=\tilde{F}_{n-2}H_{k}(\tilde{X}_Y)=\Z$, so
the result is true.

For $k=2n-1, 2n-2$, $G_{n-2}H_k(\tilde{X}_Y,
{\mathbb{Q}})=\tilde{F}_{n-2}H_k(\tilde{X}_Y,{\mathbb{Q}})=H_k(\tilde{X}_Y,
{\mathbb{Q}})$ follows from the definitions of the geometric
filtration and the Hodge filtration.

The only case left is $k=2n-3$ since the case that $k=2n-4$ has been
proved in Theorem 1.1. In this case, $T_{n-2}H_k(M,
{\mathbb{Q}})=G_{n-2}H_k(M,{\mathbb{Q}})$ has been proved in
\cite{author2} for any smooth projective variety $M$. The assumption
$G_{n-2}H_k(X,{\mathbb{Q}})=\tilde{F}_{n-2}H_k(X,{\mathbb{Q}})$ is
equivalent to $T_{n-2}H_k(X,{\mathbb{Q}}) =
\tilde{F}_{n-2}H_k(X,{\mathbb{Q}})$ in this situation. Hence
$T_{n-2}H_k(\tilde{X}_Y, {\mathbb{Q}}) =
\tilde{F}_{n-2}H_k(\tilde{X}_Y,{\mathbb{Q}})$ follows from
Proposition \ref{sec:tophdgcod2}. Now by (\ref{eq:TopGeoHodge}), we
have $G_{n-2}H_k(\tilde{X}_Y, {\mathbb{Q}}) =
\tilde{F}_{n-2}H_k(\tilde{X}_Y,{\mathbb{Q}})$.

\medskip
On the other hand, it has been proved in [\cite{Lewis1}, Lemma 13.6]
that $G_{n-2}H_k(X,{\mathbb{Q}})\cong
\tilde{F}_{n-2}H_k(X,{\mathbb{Q}})$ holds if
$G_{n-2}H_k(\tilde{X}_Y,{\mathbb{Q}})\cong
\tilde{F}_{n-2}H_k(\tilde{X}_Y,{\mathbb{Q}})$. The last part is
exactly the assumption. This completes the proof of Theorem 1.2 for
one blowup over a smooth subvariety of codimension at least two.

 \qe


\ss{The Proof of Theorem 1.3 for a Blowup}  Similarly, for 1-cycles,
we have the following.

{\Prop  \label{sec:tophdgdim1} For integer $k\geq 2$, the assertion
that ``$T_1H_k(X,{ \mathbb{Q}})=\tilde{F}_1H_k(X,{\mathbb{Q}})$
holds" is a birationally invariant property of smooth
$n$-dimensional varieties $X$.}

\medskip
\bp As before, there are two cases to consider:

\textbf{Case a:} If $T_1H_k(X,{\mathbb{Q}})=\tilde{F}_1H_k(X,{
\mathbb{Q}})$ holds, then $T_1H_k(\tilde{X}_Y, {\mathbb{Q}}) =
\tilde{F}_1H_k(\tilde{X}_Y,{\mathbb{Q}})$ holds. By the theorems in
[\cite{Friedlander-Mazur},\S 7], $T_1H_k(M,{\mathbb{Q}}) \subseteq
\tilde{F}_1H_k(M,{\mathbb{Q}})$ holds for any smooth variety $M$. We
only need to show $T_1H_k(\tilde{X}_Y,{\mathbb{Q}}) \supseteq
\tilde{F}_1H_k(\tilde{X}_Y, {\mathbb{Q}})$. The argument is similar
to the proof of the Theorem 1.3 in \cite{author2}. I give the detail
as follows:

Let $a\in \tilde{F}_{1}H_{k}(\tilde{X}_Y,{\mathbb{Q}})$, set
$b=\sigma_*(a)\in \tilde{F}_{1}H_{k}(X,{\mathbb{Q}})$. By
assumption, there exists $\tilde{b}\in L_{1}H_{k}(X,{\mathbb{Q}})$
such that $\Phi_{1,k}(\tilde{b})=b$. By the blowup formula in Lawson
homology (see \cite{author1}), we know that $\sigma_*:
L_{1}H_{k}(\tilde{X}_Y,{\mathbb{Q}})\rightarrow L_{1}H_{k}(X,{
\mathbb{Q}})$ is surjective, there exists an element $\tilde{a}\in
L_{1}H_{k}(\tilde{X}_Y,{\mathbb{Q}})$ such that
$\sigma_*(\tilde{a})=\tilde{b}$. By the commutative diagram in Lemma
\ref{eq:bpdiaglh} and Corollary \ref{sec:bpdiaghom}, we have
$j^*(\Phi_{1,k}(\tilde{a})-a)=0\in H_k^{BM}(U,{\mathbb{Q}})$. The
exactness of the localization sequence in the rows in Corollary
\ref{sec:bpdiaghom} implies that there exists an element $c\in
H_k(D,{\mathbb{Q}})$ such that $i_*(c)=\Phi_{1,k}(\tilde{a})-a$. Set
$d=\pi_*(c)\in L_{1}H_{k}(Y)\otimes{\mathbb{Q}}$. By the commutative
diagram in Corollary \ref{sec:bpdiaghom}, $d$ maps to zero under
$(i_0)_*: H_k(Y, { \mathbb{Q}})\rightarrow H_k(X,{\mathbb{Q}})$.
Hence there exists an element $e\in H_{k+1}^{BM}(U,{\mathbb{Q}})$
such that $(\delta_0)_*(e)=d$. Let $\tilde{d}\in
H_k(D,{\mathbb{Q}})$ be the image of $e$ under this boundary map
$\delta_*: H_{k+1}^{BM}(U,{\mathbb{Q}})\rightarrow
H_k(D,{\mathbb{Q}})$, i.e., $\tilde{d}=\delta_*(e)$. Therefore, the
image of $c-\tilde{d}$ is zero in $H_k(Y,{\mathbb{Q}})$ under
$\pi_*$ and is zero in $H_k(\tilde{X}_Y,{\mathbb{Q}})$ under $i_*$.
By the blowup formula in Lawson homology (see \cite{author1}), we
know such an element $c-\tilde{d}$ in the image of some $f\in
L_1H_k(D)\otimes {\mathbb{Q}}$, i.e.,$\Phi_{1,k}(f)=c-\tilde{d}$.
Hence we get $\Phi_{1,k}(\tilde{a}-i_*(f))=a$. This is the
surjectivity as we want.

\medskip
\textbf{Case b:}
 If $T_1H_k(\tilde{X}_Y, {
\mathbb{Q}})=\tilde{F}_1H_k(\tilde{X}_Y,{\mathbb{Q}})$ holds, then
$T_1H_k(X,{\mathbb{Q}})=\tilde{F}_1H_k(X,{\mathbb{Q}})$ holds. This
part is relatively easy. As before, we only need to show
$T_1H_k(X,{\mathbb{Q}})\supseteq\tilde{F}_1H_k(X,{\mathbb{Q}})$.

Let $b\in\tilde{F}_1H_k(X,{\mathbb{Q}})$. Since
$\sigma:\tilde{X}_Y\rightarrow   X $ is the blowup along the smooth
variety $Y$, we have
$\sigma_*(\tilde{F}_1H_k(\tilde{X}_Y,{\mathbb{Q}}))    \subseteq
\tilde{F}_1H_k(X,{\mathbb{Q}})$. In fact, the inclusion is an
equality. (See \cite{Lewis2} Lemma.13.6) Therefore, there is an
element $a\in \tilde{F}_1H_k(\tilde{X}_Y,{\mathbb{Q}})$ such that
$\sigma_*(a)=b$. By assumption, there is an element $\tilde{a}\in
L_1H_{k}(\tilde{X}_Y,{\mathbb{Q}})$ such that
$\Phi_{1,k}(\tilde{a})=a$. Set $\tilde{b}=\Phi_{1,k}(\tilde{a})\in
L_1H_k(X,{\mathbb{Q}}) $. By the naturality of $\Phi_{1,k}$, we have
$\sigma_*(\tilde{b})=b$. This is the surjectivity as we need.

\qe

\medskip
Now we give the proof of Theorem 1.3. First, suppose
$G_{1}H_k(X,{\mathbb{Q}})=\tilde{F}_{1}H_k(X,{\mathbb{Q}})$. We want
to show that $G_{1}H_k(\tilde{X}_Y, {
\mathbb{Q}})=\tilde{F}_{1}H_k(\tilde{X}_Y,{\mathbb{Q}})$.

Now comparing the blowup formula for Lawson homology (cf.
\cite{author1}) and for singular homology (both with ${\mathbb{Q}}$
coefficients) along the same smooth subvariety $Y$ of codimension at
least two, we find the same new components, i.e.,

$$
\bigoplus_{j=1}^{r-1}H_{k-2j}(Y,{\mathbb{Q}}),
$$
both in $L_{1}H_k(\tilde{X}_Y,{\mathbb{Q}})$ and
$H_k(\tilde{X}_Y,{\mathbb{Q}})$.

This, together with (\ref{eq:TopGeoHodge}), implies that the new
component of this blowup along $Y$ in
$G_{1}H_k(\tilde{X}_Y,{\mathbb{Q}})$ contains
$\bigoplus_{j=1}^{r-1}H_{k-2j}(Y,{\mathbb{Q}})$. Since
$G_{1}H_k(\tilde{X}_Y,{\mathbb{Q}})\subseteq H_k(\tilde{X}_Y,{
\mathbb{Q}})$, the new component of this blowup along $Y$ in
$G_{1}H_k(\tilde{X}_Y,{\mathbb{Q}})$ is also contained in
$\bigoplus_{j=1}^{r-1}H_{k-2j}(Y,{\mathbb{Q}})$. Therefore

\begin{equation}\label{eq:geoblowup}
G_{1}H_k(\tilde{X}_Y,{\mathbb{Q}})\cong
\bigg\{\bigoplus_{j=1}^{r-1}H_{k-2j}(Y,{\mathbb{Q}})\bigg\}\bigoplus
G_{1}H_k(X,{\mathbb{Q}})
\end{equation}

Similarly,

\begin{equation}\label{eq:Hodgeblowup}
\tilde{F}_{1}H_k(\tilde{X}_Y,{\mathbb{Q}})\cong
\bigg\{\bigoplus_{j=1}^{r-1}H_{k-2j}(Y,{\mathbb{Q}})\bigg\}\bigoplus
\tilde{F}_{1}H_k(X,{\mathbb{Q}})
\end{equation}

From (\ref{eq:geoblowup}) and (\ref{eq:Hodgeblowup}), we deduce that
$G_{1}H_k(\tilde{X}_Y,{
\mathbb{Q}})=\tilde{F}_{1}H_k(\tilde{X}_Y,{\mathbb{Q}})$.

\medskip
On the other hand, we also need to show that if
$G_{1}H_k(\tilde{X}_Y,{\mathbb{Q}})=\tilde{F}_{1}H_k(\tilde{X}_Y,{
\mathbb{Q}})$, then
$G_{1}H_k(X,{\mathbb{Q}})=\tilde{F}_{1}H_k(X,{\mathbb{Q}})$. An
argument similar to the one given in \textbf{Case B} works. Lewis
[\cite{Lewis1}, Lemma 13.6] proved this part in a more general
setting.

This finishes the proof of Theorem 1.3 for a blowup along a smooth
subvariety with codimension at least two.

\qe

\medskip

Now recall the weak factorization Theorem proved in \cite{AKMW} (and
also \cite{Wlodarczyk}) as follows:

{\Th (\cite{AKMW} Theorem 0.1.1, \cite{Wlodarczyk}) Let $f \colon X
\to X'$ be a birational map of smooth complete varieties over an
algebraically closed field of characteristic zero, which is an
isomorphism over an open set $U$. Then $f$ can be factored as

$$X = X_0
\stackrel{\varphi_1}{\rightarrow} X_1
\stackrel{\varphi_2}{\rightarrow}\dots
\stackrel{\varphi_{n+1}}{\rightarrow} X_n = X'$$ where each $X_i$
is a smooth complete variety, and $\varphi_{i+1}: X_i \to X_{i+1}$
is either a blowing-up or a blowing-down of a smooth subvariety
disjoint from $U$.

Moreover, if $X - U$ and $X' - U$ are simple normal crossings
divisors, then the same is true for each $X_i - U$, and the center
of the blowing-up has normal crossings with each $X_i - U$.}

\medskip
Hence {${\rm Hodge}^{2,2}(X,\mathbb{Q})$} ,
{$\widetilde{GHC}(n-2,k,X)$} and $\widetilde{GHC}(1,k,X)$ are
birationally invariant properties about the smooth manifold $X$.

 \qe

\medskip
The proof of the Corollary 1.1  and 1.2 are based on Theorem 1.1,
Remark 1.1 and the strong Lefschetz Theorem. By using the strong
Lefschetz Theorem, one can show that ${\rm Hodge}^{p,p}(X,{
\mathbb{Q}})\Rightarrow {\rm Hodge}^{n-p,n-p}(X,\mathbb{Q})$ for
$2p\leq n$. (See \cite{Lewis1} for the details.)

\qe

\medskip
The Corollary 1.3 is obvious from Theorem 1.2 and Theorem 1.3.

\qe

\medskip
\s {A Remark on Generalizations}

\medskip
From the proof of the Theorem 1.1 and 1.2, we can draw the following
conclusions:

\begin{enumerate}
\item[(a)] Fix $n>0$ and $0\leq p\leq n$. If we have
${\rm Hodge}^{i,i}(Y,\mathbb{Q})$ for all $i\leq p$ and all smooth
projective variety $Y$, i.e., the Hodge conjecture is true for every
smooth projective variety $Y$ with $\dim(Y)=n$ and for algebraic
cycles with codimension $\leq p$, then  ${\rm Hodge}^{p+1,p+1}(X,
\mathbb{Q})$ is a birational invariant statement for every smooth
projective $X$ with $\dim(X)\leq n+2$. For example, if we have ${\rm
Hodge}^{2,2}(Y,\mathbb{Q})$ for all 4-folds $Y$, then ${\rm
Hodge}^{p,p}(X,\mathbb{Q})$ is a birational statement for any
integer $0\leq p\leq \dim(X)$ and smooth projective varieties $X$
with $\dim(X)\leq 7$.

\medskip
For the Generalized Hodge Conjecture, we have
\item[(b)] Fix $n>0$ and $0\leq p\leq n$. If we have
$\widetilde{GHC}(i,k,Y)$ for $i\leq p$, i.e., the Generalized Hodge
Conjecture is true for every smooth projective $Y$ with $\dim(Y)=n$
and for algebraic cycles with codimension $\leq p$, then
$\widetilde{GHC}(m-p-1,k,X)$  is a birational invariant statement
for every smooth projective variety $X$ with $\dim(X)=m\leq n+2$.

\medskip
Similarly,
\item[(c)] Fix $n>0$ and $0\leq p\leq n$. If we have
$\widetilde{GHC}(i,k,Y)$ for $i\leq p$, i.e., the Generalized Hodge
Conjecture is true for every smooth projective $Y$ with $\dim(Y)=n$
and for algebraic cycles with dimension $\leq p$, then
$\widetilde{GHC}(p+1,k,X)$  is a birational invariant statement for
every smooth projective variety $X$ with $\dim(X)=m\leq n+2$.

As a corollary of part  (b) and (c), we have, for example, if we
have $\widetilde{GHC}(1,3,Y)$ for all 3-folds $Y$, then
$\widetilde{GHC}(p,k,X)$ is a birational statement for $X$ with
$\dim(X)\leq 5$.

\end{enumerate}

\medskip
{\Large\begin{center}{\bf Acknowledge}\end {center}}{\hskip .2 in} I
would like to express my gratitude to my advisor, Blaine Lawson, for
all his help.

\medskip

\end{document}